% Logic Eprints
%Submitted 2005 Mon Jul 11, 1994 by: nanovx!bluejay!moth (thomas e. leathrum )
%logic/leathrum/offbranch.tex
%

\input amstex
\documentstyle{amsppt}
% \magnification=\magstep1
\NoRunningHeads

%%%  ALIASES
\redefine\qed{\null\hfill$\blacksquare$\hfilneg\null\par}
\define\varqed{\null\hfill//\hfilneg\null\par}
\define\Union{\operatornamewithlimits{\bigcup}}
\define\Intersection{\operatornamewithlimits{\bigcap}}
\define\forces{\Vdash}
\define\satisfies{\vDash}
\define\restrict{\!\upharpoonright\!}
\define\concat{\!{}^\frown\!}
\define\seq#1{\left<{#1}\right>}
\define\bsl{\!\smallsetminus\!}

\redefine\epsilon{\varepsilon}
\define\ldbrack{[\kern-1.5pt[}
\define\rdbrack{]\kern-1.5pt]}
\define\decides{\,\Vert\,}
\define\modequiv{/\!\!\!\equiv\!\!\null}

%%% SPECIAL MACROS
% Macro Names:
%     \secnum
%     \Subsecnum
%     \SSubsecnum
%     \SSSubsecnum

% declare and set the counters
\newcount\seccount
\seccount0
\newcount\Subseccount
\Subseccount0
\newcount\SSubseccount
\SSubseccount0
\newcount\SSSubseccount
\SSSubseccount0
\newcount\figcount
\figcount0

% declare routines to advance and reset counters

\def\nextsec
     {\advance\seccount1
      \global\Subseccount0
      \global\SSubseccount0
      \global\SSSubseccount0
      \global\figcount0}
\def\nextSubsec
     {\advance\Subseccount1
      \global\SSubseccount0
      \global\SSSubseccount0}
\def\nextSSubsec
     {\advance\SSubseccount1
      \global\SSSubseccount0}

% declare routines to print numbers in appropriate format
\def\secnum{\global\nextsec
            {\bf\number\seccount\ \ }}
\def\Subsecnum{\global\nextSubsec
     {\rm\number\seccount.\number\Subseccount\ \ }}
\def\SSubsecnum{\global\nextSSubsec
     {\rm\number\seccount.\number\Subseccount
      .\number\SSubseccount
      \ \ }}  

\def\SSSubsecnum{\global\advance\SSSubseccount1
     {\rm\number\seccount.\number\Subseccount
      .\number\SSubseccount.\number\SSSubseccount
      \ \ }}
\def\fignum{\global\advance\figcount1{\it  
\number\seccount.\number\figcount}}

%%% END PREAMBLE

\topmatter
\title A Special Class of Almost Disjoint Families
\endtitle
\author  Thomas E.\ Leathrum
\endauthor
\address
Thomas E.\ Leathrum\newline
Department of Mathematics and Computer Science\newline
Dartmouth College\newline
Hanover, NH\qquad 03755
\endaddress
\curraddr \newline
Division of Mathematical and Natural Sciences\newline
Berry College\newline
Mount Berry, GA\qquad 30165
\endcurraddr
\email moth\@bluejay.atl.ga.us
\endemail
\abstract
The collection of branches (maximal linearly ordered sets of nodes)  
of
the tree ${}^{<\omega}\omega$ (ordered by inclusion) forms an almost
disjoint family (of sets of nodes).  This family is not maximal --  
for
example, any level of the tree is almost disjoint from all of the
branches.  How many sets must be added to the family of branches to  
make
it maximal?  This question leads to a series of definitions and  
results:
a set of nodes is {\it off-branch} if it is almost disjoint from  
every
branch in the tree; an {\it off-branch family} is an almost disjoint
family of off-branch sets; 

${\frak o}$
is the minimum cardinality of a maximal off-branch family.  

Results concerning $\frak o$ include: (in ZFC) 

${\frak a}\leq{\frak o}$, and (consistent with ZFC) $\frak o$ is not
equal to any of the standard small cardinal invariants $\frak b$, 

$\frak a$, $\frak d$, or ${\frak c}=2^\omega$.  Most of these
consistency results use standard forcing notions -- for example,
${\frak b}={\frak a}<{\frak o}={\frak d}={\frak c}$ 

in the Cohen model.  

Many interesting open questions remain, though -- for example,
$\frak d\leq\frak o$.
\endabstract
\endtopmatter

\document

%%\openup4\jot

The results in this paper have arisen from 

a study of structural
and combinatorial properties of almost disjoint families, in
particular the effects of various kinds of forcing on such families.
It is known, for example, that if $V\satisfies CH$ and $\Bbb P$ is
constructed by a finite support product of Cohen forcing, then there  
is a maximal
almost disjoint family $\Cal A$ in $V$ which remains maximal in the  
extension 

$V^{\Bbb P}$ \cite{Ku}.  Similar results can be shown with different  
assumptions --
e.g.\ if $V\satisfies MA+\neg CH$, or if $\Bbb P$ adds random reals  
instead.  On the
other hand, the collection of branches of  the tree  
${}^{<\omega}\omega$ 

form an almost disjoint family of size continuum.  This
family is not maximal, but can be easily extended using Zorn's Lemma.   
However, any
time a forcing extension adds a new real, a new branch through the  
tree is added --
and so, in the extension, the almost disjoint family is no longer  
maximal.  Two
general questions arise from such examples:  What properties of an  
almost
disjoint family in the ground model can be used to make this  
distinction more
precise?  When extending a 

particular nonmaximal almost disjoint family to a maximal family,
how many new sets must
be added?  This paper looks closely at the second question, in the 

special case given above
(the nonmaximal family being the branches through the tree  
${}^{<\omega}\omega$).

\subhead
\secnum
Basic Invariants
\endsubhead

Definitions of small cardinal
invariants have the general form \cite{vD}: 

$$
\min\{|{\Cal Q}|: {\Cal Q}\subseteq [\omega]^\omega
\ \text{is a family satisfying property $Q$}\}.
$$
This section is devoted to devising a basic list of small cardinal  
invariants, 

so it will
be necessary to define several different properties $Q$.  

For example:

\definition{\Subsecnum Definition}
Two infinite sets $A,B\subseteq\omega$ are {\it almost disjoint}
if their interesection is finite.
An {\it almost disjoint family} is a collection of
infinite subsets of $\omega$ which are pairwise almost disjoint. 

\varqed\enddefinition

Let $Q$ be the property that $\Cal Q$ is an infinite maximal almost  
disjoint family.  

Then the resulting small cardinal
invariant is called $\frak a$.

The columns of $\omega\times\omega$ form a
decomposition of $\omega\times\omega$ --- so that if $Q$ is the  
property of being
pairwise almost disjoint subsets of $\omega\times\omega$ which are  
also almost
disjoint from every column, and being maximal under this property,  
then the resulting
small cardinal invariant is still equal to $\frak a$.  So tighten up  
the property
somewhat:  require further than every set in the family intersect any  
given column
{\it at most once} (so that sets in the family can be regarded as  
infinite partial
functions from $\omega$ to $\omega$, which are pairwise almost  
disjoint as sets of
pairs).  With this stronger $Q$, the resulting small cardinal  
invariant is called
${\frak a}_s$.  It is clear now that ${\frak a}\leq{\frak a}_s$.

Some cardinal invariants are defined in terms of ${}^\omega\omega$,  
the space of
functions from the natural numbers to the natural numbers, ordered by  
``eventual
domination" $<^*$ --- $f<^* g$ if and only if 

$\{n:g(n)\leq f(n)\}$ is finite.  For example, if
the property $Q$ over  ${}^\omega\omega$ is the property of being  
unbounded in this
ordering, then the resulting small cardinal invariant is called  
$\frak b$.

\definition{\Subsecnum Definition}
 A family ${\Cal D}\subseteq {}^\omega\omega$ is 

{\it dominating} if for every $f\in{}^\omega\omega$ there is a  
$g\in{\Cal D}$ such that
$f<^* g$.
\varqed\enddefinition

If the property $Q$ is the property of being dominating, then the  
resulting
small cardinal invariant is called $\frak d$.  Since any dominating  
family is unbounded,
it is clear that ${\frak b}\leq{\frak d}$.

There are many inequalities provable between these cardinal  
invariants --- for example,
the inequalities ${\frak a}\leq{\frak a}_s$ 

and ${\frak b}\leq{\frak d}$ have already been
mentioned.  

The only other inequality known in $ZFC$ for these invariants is
$\frak b\leq\frak a$ \cite{vD}.
There are many other invariants which have been investigated, and 

many inequalities between them have been established.
  Some of these results are quite difficult --- the interested reader  
is
referred to \cite{vD} or \cite {BS} for details.

By forcing techniques, it is possible to construct models in which  
various {\it strict}
inequalities hold between these cardinal invariants.  For example,  
forcing to add
$\omega_2$-many Cohen reals over a model of $CH$ gives a model of 

${\frak a}_s<{\frak d}$.  Some such consistency questions remain  
open, though 

--- for example, it is unknown whether $\frak a\leq\frak d$ is  
provable
in $ZFC$.  In fact, this question is an important motivation for the
present research.

\subhead\secnum Definitions
\endsubhead

\definition{\Subsecnum Definition}
The set ${}^{<\omega}\omega$ consists of finite sequences of natural  
numbers.  This set
is given an ordering by $\sigma\leq \tau$ if and only if $\sigma$ is  
an initial
segment of $\tau$
--- the result is a countably branching, countable
height tree ordering.  A {\it node} of the tree is an element of  
${}^{<\omega}\omega$.
A {\it branch} through the tree is a maximal linearly ordered
set of nodes.  The {\it $n^{th}$ level} of the tree is the set of  
nodes which, as
sequences, all have length $n$.
\varqed\enddefinition

In general, the families this paper deals with will be almost  
disjoint families of
infinite sets of nodes, in particular families extending the
(nonmaximal) family of branches of ${}^{<\omega}\omega$.

\definition{\Subsecnum Definition}
An infinite set $A$ of nodes of ${}^{<\omega}\omega$ is {\it  
off-branch} if $A$ is
almost disjoint from every branch of ${}^{<\omega}\omega$.  An {\it  
off-branch family}
is an almost disjoint family of off-branch sets.
\varqed\enddefinition

\definition{\Subsecnum Definition}
$$
{\frak o}=\min\{|{\Cal O}|:{\Cal O}\ \text{is a maximal off-branch  
family}\}.
$$
\varqed\enddefinition

\subhead\secnum Equivalent and Related Invariants
\endsubhead

The first few results concern equivalent definitions of $\frak o$.   
For example, one
natural question is whether it makes any difference to define $\frak  
o$ in terms of the
binary tree ${}^{<\omega}2$ instead of ${}^{<\omega}\omega$.  In  
order to establish
the context for this question, define a set $A$ of nodes of the tree 

${}^{<\omega}2$ to be {\it off-binary} if $A$ is almost disjoint from  
every branch of 

${}^{<\omega}2$, and an {\it off-binary family} is an almost disjoint  
family of
off-binary sets.  By analogy with the definition of $\frak o$, define
$$
\frak o_b=\min\{|{\Cal O}|:  {\Cal O}\ \text{is a maximal off-binary  
family}\}.
$$

\proclaim{\Subsecnum Lemma}  $\frak o=\frak o_b$.
\endproclaim

\demo{Proof}
The basic idea of this proof is to construct mappings between  
${}^{<\omega}2$ and 

${}^{<\omega}\omega$ which carry off-binary families to off-branch  
families, and vice
versa.  It turns out that the best thing to do is to simply embed the  
trees in
canonical ways into each other, and look at {\it pullbacks} of the  
families.

For one direction, notice first that  
${}^{<\omega}2\subseteq{}^{<\omega}\omega$, so
the identity map on ${}^{<\omega}2$ embeds it into  
${}^{<\omega}\omega$.  Let $\Cal O$
be a maximal off-branch family.  

Let 

$$
\align
\bar\Cal O=\{A\cap{}^{<\omega}2:&  A\in\Cal O\ \text{and}\\
& A\cap{}^{<\omega}2\ \text{infinite}\}.
\endalign
$$
(So $\bar\Cal O$ is the pullback of $\Cal O$ over the identity map  
embedding.)
Since each $A\in\Cal O$ is off-branch, $\bar\Cal O$ is an off-binary  
family.  If
$\bar\Cal O$ is not maximal, let $B\subseteq{}^{<\omega}2$ be a  
witness to this fact
--- i.e.\ an off-binary set which is almost disjoint from every  
element of 

$\bar\Cal O$.  Then $B\subseteq{}^{<\omega}\omega$, $B$ is an  
off-branch set, and $B$
is almost disjoint from every element of $\Cal O$, contradicting  
maximality of 

$\Cal O$.  This proves, in particular, that $\frak o_b\leq \frak o$.

To do the other direction, define an embedding 

$\pi:{}^{<\omega}\omega\to{}^{<\omega}2$ as follows.  

Let ${\vec 1}_n$ be the length $n$ sequence of 1's.  For 

$\sigma\in{}^{<\omega}\omega$, $\sigma=\seq{n_0,n_1,\dots,n_k}$, let
$$
\pi(\sigma)={\vec 1}_{n_0}\concat\seq{0}\concat
{\vec 1}_{n_1}\concat\seq{0}\concat{}\dots{}\concat
{\vec 1}_{n_k}\concat\seq{0}.
$$
So the image of $\pi$ is the collection of binary sequences ending in  
0, with the
coordinates of $\sigma$ being coded by the lengths of corresponding  
blocks of
consecutive 1's.

With this $\pi$, the argument proceeds much as in the first case.   
Let $\bar\Cal O$ be
a maximal off-binary family.  For $\bar A \in \bar\Cal O$, let 

$A=\{\sigma\in{}^{<\omega}\omega:\pi(\sigma)\in\bar A\}$ (so $A$ is  
the pullback
of $\bar A$ over $\pi$), and let 

$\Cal O=\{A:\bar A\in\bar\Cal O, A\ \text{infinite}\}$.  Then $\Cal  
O$ is
an off-branch family.  Furthermore, if $\Cal O$ is not maximal and 

$B\subseteq{}^{<\omega}\omega$ witnesses the fact that $\Cal O$ is  
not maximal, then
$\pi''B=\{\pi(\sigma):\sigma\in B\}$ 

witnesses that $\bar\Cal O$ isn't maximal either, which contradicts  
the
assumption.  So $\frak o\leq \frak o_b$.
\qed\enddemo

Another equivalent form of $\frak o$ adds the requirement that the  
off-branch family
contain a decomposition of the tree.
Define by analogy with $\frak o$ the cardinal invariant
$$
\align
\frak o_d=\min\{|\Cal O|:\ &\text{ there is a decomposition $\Cal D$  
of 
${}^{<\omega}\omega$ into infinite sets so that}\\
&\Cal O\cup\Cal D
\ \text{is a maximal off-branch family}\}.
\endalign
$$

\proclaim{\Subsecnum Lemma}  $\frak o=\frak o_d$.
\endproclaim

\demo{Proof}
Clearly $\frak o\leq \frak o_d+\omega = \frak o_d$.  So it only  
remains to show that 

$\frak o_d\leq \frak o$.
Let $f:\omega\to\Cal O$ be an injection, and let  
$g:\omega\to{}^{<\omega}\omega$
be a bijection.  Define  a function $h$ by 

$$
h(n)=(\{g(n)\}\cup f(n))\bsl\Union_{k<n} f(k).
$$
Then the range of $h$ is a decomposition of ${}^{<\omega}\omega$ into  
infinite 

off-branch sets, and for any off-branch set $A$, $h(n)\cap A$ is  
infinite if and only
if $f(n)\cap A$ is infinite.  Letting $\Cal O'=\Cal O\bsl ran(f)$  
provides a witness
to $\frak o_d\leq \frak o$.
\qed\enddemo

The situation changes somewhat when considering a {\it particular}
decomposition of the tree.
Given a decomposition $\Cal D$ of  of 

${}^{<\omega}\omega$ into  off-branch sets,
define, by analogy with $\frak o$, the cardinal invariant
$$
\tilde\frak o(\Cal D) = \min\{|\Cal O|: \Cal O\cup \Cal D
\ \text{is a maximal off-branch family}\}.
$$
So $\frak o_d=\min_{\Cal D} \frak o(\Cal D)$.
Clearly $\frak o_d\leq\frak o(\Cal D)$ for any $\Cal D$ --- 

but what about $\frak o(\Cal D)\leq\frak o_d$? 

It turns out that it is easier to approach this problem directly, in  
terms of 

$\frak o$ rather than $\frak o_d$.  However, this result requires the  
use of the fact
that $\frak a\leq \frak o$, which will be proved later.

\proclaim{\Subsecnum Lemma} 

For any decomposition $\Cal D$ of  of 

${}^{<\omega}\omega$ into  off-branch sets,
$\frak o(\Cal D)=\frak o$.
\endproclaim

\demo{Proof} Deferred to next section.
\enddemo

Antichains in the tree are clearly off-branch --- they intersect any  
given branch at
most once.  So what happens to the cardinal invariant if you  
strengthen the off-branch
condition to talking about antichains?  (Notice that every infinite  
off-branch set
contains an infinite antichain.)
Define, by analogy with $\frak o$, yet another
cardinal invariant:
$$
\bar\frak o = \min\{|\Cal O|: \Cal O\ \text{is a maximal almost  
disjoint family of
antichains of ${}^{<\omega}\omega$}\}.
$$
Again, clearly $\frak o\leq\bar\frak o$.  While it seems  
counterintuitive that
$\bar\frak o\leq\frak o$, constructing a model in which these two  
invariants are
different seems quite difficult.

\subhead\secnum Results in $ZFC$ --- $\frak a \leq \frak o$
\endsubhead

Recall that $\frak a$ is the minimum size of a maximal almost  
disjoint family of
subsets of a countable set.  It is not obvious at first that $\frak  
o$ is related to
$\frak a$, since the off-branch family does {\it not} contain the  
branches.

\proclaim{\Subsecnum Theorem}  $\frak a \leq \frak o$.
\endproclaim

(The following version of this proof was suggested by the referee to  
show
more clearly the connection with $\frak a_s$ later.)

\demo{Proof}
Suppose, for contradiction, that $\frak o<\frak a$.
Let $\Cal O$ be a maximal off-branch family of size $\frak o$.
For each $n$, let $b_n$ be any branch containing the node $\seq{n}$.  

Then the $b_n$
are all distinct, and in fact disagree at the very first level.
For each $O\in\Cal O$, define:
$$
\bar O=\{\seq{n,i}:\exists \sigma\in O\cap b_n\ ht(\sigma)=i\}.
$$
Let $Col$ be the collection of columns of $\omega\times\omega$ ---
then $\{\bar O: O\in\Cal O,\ \bar O\ \text{infinite}\}\cup Col$
is an almost disjoint family of subsets of $\omega\times\omega$.
However, this family cannot be maximal, since it has size $<\frak a$.
Let $\bar B\subseteq\omega\times\omega$ be almost disjoint from every
column and from every $\bar O$.
Define $B=\{\sigma:\seq{n,i}\in\bar B,\ \sigma\in b_n,\  
ht(\sigma)=i\}$.
Then $B$ is off-branch and almost disjoint from every $O\in\Cal O$,
contradicting maximality of $\Cal O$.
\qed\enddemo

It is now possible to prove Lemma 3.3, that $\frak o=\frak o(\Cal  
D)$.

\demo{\Subsecnum Proof of Lemma 3.3}
Since $\frak o=\frak o_d\leq\frak o(\Cal D)$,  it only remains to  
show
$\frak o(\Cal D)\leq\frak o$.
Let $\Cal O$ be a maximal off-branch family of size $\frak o$, 

let $\Cal D=\{D_n:n<\omega\}$ be a decomposition of  
${}^{<\omega}\omega$
into off-branch sets,
and for each $A\in\Cal O$ define
$\Cal D\restrict A=\{D_n\cap A:n<\omega, D_n\cap A\  
\text{infinite}\}$.  

Then $\Cal D\restrict A$ is a countable
pairwise disjoint collection of subsets of $A$.  So let $\Cal B_A$ be  
such that
$(\Cal D\restrict A)\cup \Cal B_A$ is a maximal almost disjoint  
family {\it of subsets
of $A$}, and $|\Cal B_A|=\frak a$.  

Since each $A$ is off-branch, each element of $\Cal B_A$ is also  
off-branch,
so $\Cal B=\Union_{A\in\Cal O} \Cal B_A$ is an off-branch family.
Also, for each $B\in\Cal B$ and each $n<\omega$, $B\cap D_n$ is  
finite,
so $\Cal B\cup\Cal D$ is an off-branch family.
Furthermore, $|\Cal B|=\frak o\cdot\frak a=\frak o$
(since $\frak a\leq\frak o$), so it remains only to show that
$\Cal B\cup\Cal D$ is a {\it maximal} off-branch family.

To show this,  let $C$ 

be an off-branch set.  

Since $\Cal O$ is a maximal
off-branch family, there is some $A\in\Cal O$ such that $C\cap A$ is  
infinite.
But $(\Cal D\restrict A)\cup \Cal B_A$ is also maximal,
so $C$ has infinite intersection with some element of 

$(\Cal D\restrict A)\cup \Cal B_A$, and thus with some element of
$\Cal B\cup\Cal D$.
\qed\enddemo

Consider now the relationship between $\frak o$ and ${\frak a}_s$.
The important distinction is
not in the nature of the decomposition, but in how often a set  
intersects one of the
specified sets --- ${\frak a}_s$
(where the sets intersect columns at most once)
bears a much closer relationship to $\bar\frak o$ (where the sets  
intersect branches
at most once) than to $\frak o(\Cal D)$.

\proclaim{\Subsecnum Theorem}  ${\frak a}_s\leq\bar\frak o$.
\endproclaim

\demo{Proof}
This proof is exactly like the proof that $\frak a\leq\frak o$, with  
the
following exception:  Notice that (using the same notation as above)
$\bar O$ now intersects each column of $\omega\times\omega$ {\it at  
most 

once}.  Therefore, the family 

$\{\bar O: O\in\Cal O,\ \bar O\ \text{infinite}\}$ is a family as in
$\frak a_s$.  The remaining details are left to the reader.
\qed\enddemo

\subhead\secnum Cohen Reals --- $Con(\frak a<\frak o=\frak d)$
\endsubhead

The  Cohen forcing notion used here will be
${\Bbb P}={}^{<\omega}\omega$
(finite sequences of natural numbers), ordered by end-extension
--- in other words, the tree turned upside-down.
This is equivalent, as a forcing notion, to the usual $Fn(\omega,2)$
Cohen forcing.

\proclaim{\Subsecnum Theorem}  

If $\Cal O$ is an off-branch family (in the ground model $V$), $\Bbb  
P$ is
the forcing notion given by the set ${}^{<\omega}\omega$ ordered by
$p\leq q$ if and only if $p\supseteq q$, 

and $G$ is $\Bbb P$-generic over $V$, then in the
extension $V[G]$, $\Cal O$ is not maximal.
\endproclaim

\demo{Proof}
For a sequence $\seq{n_0,n_1,\dots,n_k}\in {}^{<\omega}\omega$,  
define
$$
rs(\seq{n_0,n_1,\dots,n_k})=\seq{n_0,n_1,\dots,n_k+1}.
$$
This is a ``right shift"
function in the tree --- if the successors of  
$\seq{n_0,n_1,\dots,n_{k-1}}$ are ordered
according to $n_k$, left to right, then $rs$ shifts  
$\seq{n_0,n_1,\dots,n_k}$ one 

place to the right.
The Cohen-generic set $G$ will simply (and exactly) be a branch  
through the tree
${}^{<\omega}\omega$.  The
idea of this proof is to show that the ``hair" of the branch $G$, 

$\bar G=\{rs(\sigma): \sigma\in G\}$, is off-branch (which is  
obvious, since it is an
antichain) and almost disjoint from all members of the off-branch  
family $\Cal O$.

This is accomplished by a density argument.  The idea is to show  
that, for
every $A\in\Cal O$, the set of sequences whose right shifts are {\it  
not} in $A$
contains a dense {\it open} set in the forcing notion (which is just  
the upside down
tree).  Formally, define
$$
D_A=\{\sigma\in{}^{<\omega}\omega: 
\forall \tau\in{}^{<\omega}\omega\ rs(\sigma\concat\tau)\notin A\}.
$$
It remains to show that $D_A$ is dense for each $A\in\Cal O$.  This  
is accomplished by
an appeal to the fact that $A$ is an off-branch set.

To see how this works, let $G$ be a generic set for the Cohen  
forcing.  So for each
$A\in\Cal O$ there is a $\sigma\in G\cap D_A$.  Now if $\tau\in G$ is  
such that
$rs(\tau)\in A$, then $\tau\subseteq\sigma$, by the definition of  
$D_A$.  In
particular, then, $\bar G \cap A$ is finite.
\qed\enddemo

\proclaim{\Subsecnum Theorem} 

If $\Bbb P$ is a forcing notion constructed as a finite support  
product of
$\kappa$-many Cohen forcing notions (where $\kappa$ is an 

uncountable cardinal) and $G$ is $\Bbb P$-generic over $V$,  

then in the extension $V[G]$, there are no
maximal off-branch families of size less than $\kappa$.
\endproclaim

The technique in this proof is a completely standard way of dealing  
with products of
Cohen reals \cite{Ku} --- nothing more will be said here.

\proclaim{\Subsecnum Corollary}  

$Con(\omega_1=\frak b=\frak a<\frak o=\frak d= 2^\omega =\kappa)$.
\endproclaim

\demo{Proof}
Start with a model $V\satisfies ZFC+CH$.  Force to add $\kappa$-many  
Cohen reals. 

By the above Theorem, together with previous results concerning the  
effects of Cohen
forcing on the other cardinal invariants \cite{vD, BS}, 

this will give the desired model.
\qed\enddemo

\subhead\secnum Random Reals --- $Con(\frak d<\frak o)$
\endsubhead

The random forcing notion used here will be
the Boolean algebra given by the measure algebra on 

the product space $2^\kappa$ --- this algebra, used as a forcing
notion, is said to {\it add $\kappa$-many random reals.}
The following standard fact for random forcing (see \cite{So} or  
\cite{Je})
will not be proven here.

\proclaim{\Subsecnum Fact}
Let $\Bbb B$ be the measure algebra on the product space $2^\kappa$,  
and let $G$ be 

$\Bbb B$-generic over $V$.  If $r\in V[G]$ is an countable subset of  
$V$,
then there is a countable set $X\subseteq\kappa$ such that 

$r\in V[G\restrict X]$.
\endproclaim

\demo{\Subsecnum Note}
The forcing notion $\Bbb B$ which adds $\kappa$ many random reals has  
$ccc$.
Since $\kappa$ is assumed to be infinite,
the measure algebra on $2^\kappa$ is equivalent (as a forcing notion)  
to the
measure algebra on $2^{\kappa\times\omega}$.
Then a generic set $G$ gives a function $f_G:\kappa\times\omega\to  
2$, and thus $\kappa$
many reals (subsets of $\omega$), $r_\alpha$ for $\alpha<\kappa$,
such that $n\in r_\alpha$ if and only if $f_G(\alpha, n)=1$.
Now the measure algebra on $2^{\omega\times\omega}$ 

adds countably many random reals, and
the measure algebra on $2^{1\times\omega}$ is said to {\it add a  
single random real} -- 

but these two measure algebras are equivalent (as forcing notions).
\enddemo

It is now possible to apply random forcing in the context of  
off-branch families.
The author owes much of the credit for  the mode of presentation
for this proof to the guidance and insistence of Prof.\ James E.\  
Baumgartner.

\proclaim{\Subsecnum Theorem}  

If $\Cal O$ is an off-branch family (in the ground model $V$), $\Bbb  
B$ is
the measure algebra on $2^\omega$ (so that forcing with this
Boolean algebra adds a single random real),
and $G$ is $\Bbb B$-generic over $V$,  then in the
extension $V[G]$, $\Cal O$ is not maximal.
\endproclaim

This proof is essentially the same as the Cohen reals proof, but the  
new context of
random reals makes some of the calculations more delicate.  

It will be carried out using off-binary families instead of
off-branch families, but these are entirely equivalent since 

$\frak o=\frak o_b$.

\demo{Proof}
The set of branches through the tree ${}^{<\omega}2$ can be  
identified with the
product space $2^\omega$, so the sets of
branches inherit a measure algebra structure from the interval.  As  
Boolean algebras,
and thus as a forcing notions, these measure algebras are equivalent.   
So rather than
forcing with the measure algebra on $2^\omega$, use the 

equivalent measure algebra on the branches
through ${}^{<\omega}2$.  Given a
sequence $\sigma\in{}^{<\omega}2$, associate with $\sigma$ the Baire  
interval,
denoted $[\sigma]$, of branches which contain $\sigma$ (all the  
branches
with $\sigma$ as a common stem), with measures for these 

Baire intervals given by
$$
m([\sigma])=
2^{-len(\sigma)}.
$$
It is easy to see that this is consistent with the measure on  
$2^\omega$.
Finally, rather than using the Boolean algebra, this proof deals  
explicitly with the
pre-order on Baire sets (elements of the $\sigma$-algebra generated  
by the Baire
intervals), where the ordering is by almost containment:
$p\leq q$ if and only if $m(p\bsl q)=0$, where $p$ and $q$ are Baire  
sets.  

Thus, countable unions and intersections of Baire sets correspond to 

countable infima and suprema (respectively) in the Boolean algebra.

As with the single Cohen real proof, the idea of this proof is to  
look at the ``hair" 

of the generic branch added by forcing with this measure algebra.  

If  $G\subseteq{\Bbb B}$ is a generic ultrafilter, 

define $g=\Union\{\sigma:[\sigma]\in G\}$ --- so $g$ is in  
${}^\omega2$, and is
essentially the generic branch added by forcing with $\Bbb B$.   
Define the ``hair" on
$g$ to be the set 

$H=\{\tau\concat\seq{1-i}:\exists n\ g\restrict n  
=\tau\concat\seq{i}\}$.  

Then $H$ is
an antichain in ${}^{<\omega}2$, so it remains only to show that $H$  
is almost
disjoint from every element of $\Cal O$.

For each $A\in{\Cal O}$, let  

$f:\omega\to A$ enumerate $A$.
First, notice that
$$
\Union_n \sum_{k\geq n} [f(k)] = \emptyset
$$
To see this,
suppose it is not true and let $b\in\Union_n \sum_{k\geq n} [f(k)]$.
Now $b\in 2^\omega$ --- interpret $b$ as an element of ${}^\omega2$,
so that $\{b\restrict i:i<\omega\}$ is a branch through the tree.
Since there are infinitely many $k$ for which $b\in[f(k)]$,
there are infinitely many $k$ for which $b\restrict i=f(k)$ for some  
$i$.
But this contradicts $A$ being off-branch.

The sets $\Union_{k\geq n} [f(k)]$ are nested (i.e.\ 

$\Union_{k\geq n} [f(k)]\supseteq\Union_{k\geq n+1} [f(k)]$ for all  
$n$),
so 

$$
m\biggl(\Intersection_n \Union_{k\geq n} [f(k)]\biggr)
=\lim_{n\to\infty} m\biggl(\Union_{k\geq n} [f(k)]\biggr) =0.
$$
For each $n$, let $B_n\subseteq\{f(k):k\geq n\}$ be the antichain of
minimal elements (with respect to the tree ordering)
of the set $\{f(k):k\geq n\}$.  Then 

$\Union_{k\geq n} [f(k)] = \Union_{\sigma\in B_n} [\sigma]$.
Since the intervals $[\sigma]$ are {\it disjoint} for $\sigma\in  
B_n$,
$$
m\biggl(\Union_{\sigma\in B_n} [\sigma]\biggr)
= \sum_{\sigma\in B_n} m([\sigma]).
$$
So
$$
\lim_{n\to\infty} \sum_{\sigma\in B_n} m([\sigma]) = 0.
$$

For all $\sigma\in{}^{<\omega}2$, define
$$
pred(\sigma)=\cases
\tau\qquad&\sigma=\tau\concat\seq{0}\ \text{or}\  
\sigma=\tau\concat\seq{1};\\
\seq{}\qquad&\text{otherwise}.
\endcases
$$
Now for all $k\geq n$, there is a $\sigma\in B_n$ such that
$[pred(f(k))]\subseteq [pred(\sigma)]$,
and in this case $m([pred(f(k))])\leq 2m([\sigma])$.  Therefore,
$$
m\biggl(\Union_{k\geq n} [pred(f(k))]\biggr)
\leq 2\cdot\sum_{\sigma\in B_n} m([\sigma])
$$
Since the sets $\Union_{k\geq n} [pred(f(k))]$ are also nested,
$$
m\biggl(\Intersection_n \Union_{k\geq n} [pred(f(k))]\biggr)
= \lim_{n\to\infty} m\biggl(\Union_{k\geq n} [pred(f(k))]\biggr)
\leq \lim_{n\to\infty} 2\cdot\sum_{\sigma\in B_n} m([\sigma]) = 0
$$
The Baire set $p=\Intersection_n \Union_{k\geq n} [pred(f(k))]$
is in the equivalence class represented by the Boolean value
$\ldbrack\forall n\ \exists k\geq n\ \exists i\ pred(f(k))
=g\restrict i\rdbrack$,
so the statement that $m(p)=0$ means 

that this value is $=\bold 0$ in the Boolean algebra.

To see why this suffices to finish the proof, notice that
$$
\forall q\ \exists n\ 
m\biggl(\Union_{k\geq n}\ [ pred(f(k)) ]\biggr)
<m(q).
$$
Let $q'=q\bsl\Union_{k\geq n}\ [ pred(f(k)) ]$.  Then for all $m$ and  
for 

all $k\geq n$, 

$q'\modequiv\forces g\restrict m \neq pred(f(k))$.  But notice that  
if 

$\sigma\in H\cap A$, then there is an $m$ such that 

$pred(\sigma)=g\restrict m$ ---
thus, since there are only finitely many $m$ such that 

$g\restrict m = pred(f(k))$ for some $k$, it
follows that $H\cap A$ is also finite.
\qed\enddemo

\proclaim{\Subsecnum Theorem}  

If $\Bbb B$ is a Boolean algebra forcing notion
which adds $\kappa$-many random reals
where $\kappa$ is an uncountable cardinal, 

and $G$ is $\Bbb B$-generic over $V$, 

then in the
extension $V[G]$, there are no maximal off-branch families of size  
less than
$\kappa$.
\endproclaim

As with Cohen reals, the technique of this proof is standard ---
see \cite{So} or \cite{Je}.

\proclaim{\Subsecnum Corollary}  

$Con(\omega_1=\frak b=\frak a=\frak d<\frak o=2^\omega=\kappa)$.
\endproclaim

\demo{Proof}
Start with a model $V\satisfies ZFC+CH$.  Force to add $\kappa$-many  
random reals. 

By the above Theorem, together with previous results concerning the  
effects of random
forcing on the other cardinal invariants \cite{vD, BS}, 

this will give the desired model.
\qed\enddemo

\subhead\secnum Sacks Reals --- $Con(\frak o<2^\omega)$
\endsubhead

\definition{\Subsecnum Definition} {\it (Sacks forcing)}
The partial ordering $\Bbb P$ is the set of perfect binary trees ---  
i.e.
$$
\align
{\Bbb P}=
\{p: &p\subseteq{}^{<\omega}2, 
\forall\sigma\in p\ \forall\tau\subseteq\sigma\ \tau\in p,  
\text{and}\\
&\forall\sigma\in p\ \exists\tau\in{}^{<\omega}2
\ \text{both}\ \sigma\concat\tau\concat\seq{0}\in p\ \text{and}
\ \sigma\concat\tau\concat\seq{1}\in p\};\\
p\leq q\ &\text{iff}\ p\subseteq q.
\endalign
$$
\varqed\enddefinition

Let $p$ be a Sacks condition (i.e.\ a perfect binary tree).  

The set $split(p)=\{\sigma\in p:\sigma\concat\seq{0}\in p\
\text{and}\ \sigma\concat\seq{1}\in p\}$ contains the ``branching  
points" or ``splitting
points" of $p$.  The set 

$split_n(p)=\{\sigma\in split(p):|\{\sigma'\subseteq\sigma:\sigma'\in  
split(p)\}|=n\}$
contains the ``$n^{th}$ branching (or splitting) points" of $p$.
For two Sacks conditions $p$ and $q$,
define $p\leq_n q$ if and only if $p\leq q$ and $split_n(q)=  
split_n(p)$ 

--- that is, $p$ agrees with $q$ up to and
including the $n^{th}$ branching points.
These orderings can be used to show that Sacks forcing satisfies 

Axiom A (see \cite{Je}) --- in particular, if $p_{n+1}\leq_n p_n$ for
each $n$ and $q=\cap_n p_n$, then $q$ is a perfect binary tree and
$q\leq_n p_n$ for each $n$.

The following result deals with antichain families instead of  
off-branch families, in
order to get a stronger result --- that forcing to add many Sacks  
reals over a model of
$CH$ leaves $\bar o=\omega_1$ in the extension.
Since this proof deals with {\it preserving} ground model families  
rather than
destroying them as in the Cohen and random reals proofs, the case of  
adding a single
Sacks real will not be handled separately --- it would simply repeat
much of the same argument as below.

\proclaim{\Subsecnum Theorem}  ($CH$)  

Suppose $\Bbb P$ is the notion of forcing which adds $\kappa$ many  
Sacks reals
with a countable support product,
where $\kappa$ is an uncountable cardinal with $cf(\kappa)>\omega$.
There is a maximal antichain family $\Cal O$ in $V$ such that $\Cal  
O$ remains
maximal in $V^{\Bbb P}$.
\endproclaim

The following technical lemmas for Sacks forcing will be useful here
(see \cite{Je}):

\proclaim{\SSubsecnum Lemma}
If $\Bbb P_\kappa$ is the forcing to add $\kappa$-many Sacks reals  
with a countable
support product, $\Bbb P_\omega$ is the forcing to add $\omega$-many  
Sacks reals (with
full support product), $G$ is $\Bbb P_\kappa$-generic over $V$, and 

$\rho\in V[G]$ is an
countable subset of $V$, then there is a $G'\in V[G]$ 

which is $\Bbb P_\omega$-generic over
$V$ such that $\rho\in V[G']$.
\endproclaim

For the purposes of the next lemma, it will be convenient to define a  
notion of
``splitting point" for $\Bbb P_\omega$.  

For $p\in\Bbb P_\omega$, an $n^{th}$ {\it splitting point} of $p$ is  
a sequence
$\seq{\sigma_0,\dots,\sigma_n}$ of elements of ${}^{<\omega}2$ such  
that
for each $i\leq n$, $\sigma_i\in split_{n-i}(p(i))$.  Then
$split_n(p)=\{\vec\sigma:\vec\sigma\ \text{an $n^{th}$ splitting  
point of $p$}\}$
and $split(p)=\Union_n split_n(p)$.
Finally $p\restrict \vec\sigma$  is defined by
$$
(p\restrict \vec\sigma)(i)=\cases
p(i)\restrict \sigma_i\quad &i\leq n;\\
p(i)\quad&\text{otherwise}.
\endcases
$$

\proclaim{\SSubsecnum Lemma}
If $p\in\Bbb P_\omega$ and $\tau$ is a $\Bbb P_\omega$-name for an  
infinite subset of
$\omega$, then there is a $q\leq p$ such that for every $n$ and every 

$\vec\sigma\in split_n(q)$, $q\restrict\vec\sigma\decides n\in\tau$.
\endproclaim

Now $q$ gives a function $\pi:split(q)\to 2$ by $\pi(\vec\sigma)=1$  
if and only if
$q\restrict\vec\sigma\forces n\in\tau$.
Such a $q$ is said to {\it code $\tau$ as $\pi$}.

\demo{Proof of Theorem}
Enumerate ${}^{<\omega}\omega$ as $\{\eta_n:n<\omega\}$.
Then the two Lemmas apply also to (names for) infinite off-branch  
sets.
By the first Lemma, it suffices to consider only $\Bbb P_\omega$ ---
if $\Cal O$ is a maximal antichain family in $V$ which is  
indestructible under
forcing with $\Bbb P_\omega$, then $\Cal O$ will also be  
indestructible under forcing
with $\Bbb P_\kappa$.  In $\Bbb P_\omega$, given names $\tau_1$ and  
$\tau_2$
(for infinite antichains)
suppose $q$ codes both $\tau_1$ and $\tau_2$ as $\pi$.
Then $q\forces\tau_1=\tau_2$.  (This determines an ``isomorphism"  
between 

$\tau_1$ and $\tau_2$.)  By the second Lemma, the set
$D=\{q:\text{$q$ codes some name $\tau$ for an infinite antichain as  
some $\pi$}\}$  

is dense in $\Bbb P_\omega$.
By $CH$, $|\Bbb P_\omega|=\omega_1$, so it is possible to enumerate  
all pairs
$\seq{(q_\alpha,\pi_\alpha):\alpha<\omega_1}$
where $q_\alpha\in D$ and $\pi:split(q_\alpha)\to 2$.  

For each $\alpha<\omega_1$, choose a name
$\tau_\alpha$ such that $q_\alpha$ codes $\tau_\alpha$ as  
$\pi_\alpha$.
(So the $\tau_\alpha$ are representative members of their  
``isomorphism classes.")
For any $p\in\Bbb P_\omega$ and any $\Bbb P_\omega$-name $\tau$ 

for an infinite antichain, there is an
$\alpha$ such that $q_\alpha\leq p$ and  
$q_\alpha\forces\tau=\tau_\alpha$.

Inductively construct antichains $A_\alpha$ for $\alpha<\omega_1$ to  
satisfy
\roster
\item $\forall\beta<\alpha\ |A_\beta\cap A_\alpha|<\omega$;
\item if
$$
q_\alpha\forces``\tau_\alpha\ \text{an infinite antichain}"
\ \text{and}
\ \forall\beta<\alpha\ q_\alpha\forces|\tau_\alpha\cap  
{A_\beta}|<\omega,
\tag{$*$}
$$
then there is some $q\leq q_\alpha$ such that
$$
q\forces|\tau_\alpha\cap A_\alpha|=\omega.
$$
\endroster
This construction establishes the context for a density argument:
Suppose the construction works, so that \therosteritem{1} and  
\therosteritem{2} hold.  

Given a name $\tau$ for an infinite antichain,
consider the set $D_\tau=\{q:\exists\xi\ q\forces \tau\cap A_\xi\  
\text{infinite}\}$.
For any $p\in\Bbb P_\omega$
let $\alpha$ be such that $q_\alpha\leq p$ and  
$q_\alpha\forces\tau=\tau_\alpha$.
Now for this pair $(q_\alpha,\tau_\alpha)$, if $(*)$ fails then there  
is an 

$\eta<\alpha$
and a $q\leq q_\alpha$
such that $q\forces\tau\cap A_\eta\ \text{infinite}$ (i.e.\ $q\in  
D_\tau$),
while if $(*)$ holds then there is a $q\leq q_\alpha$ such that 

$q\forces\tau\cap A_\alpha\ \text{infinite}$ (i.e.\ $q\in D_\tau$).
So $D_\tau$ is dense in $\Bbb P_\omega$.
The only difficulty then is in constructing $q$ when $(*)$ holds.

Enumerate $\{A_\beta:\beta<\alpha\}$ as $\{B_i:i<\omega\}$.
By ($*$), for each $i<\omega$,
$$
q_\alpha\forces |\tau_\alpha\bsl(B_0\cup\dots\cup B_i)|=\omega.
$$
The idea of this proof is to construct another ``fusion sequence"
$\seq{p_n:n<\omega}$ with $p_0=q_\alpha$ and $p_{n+1}\leq p_n$,  
constructing along 

with $p_n$
a finite ``block" $S_n\subseteq {}^{<\omega}\omega$ of $A_\alpha$.

At stage 0 of this construction, let $p_0=q_\alpha$ and  
$S_0=\emptyset$.
It will be necessary to carry several hypotheses through the  
induction --- 

in particular, after completing stage $n+1$,
$\Union_{m\leq n+1} S_m$ must be an antichain and $p_{n+1}$ must  
force that infinitely
many nodes $\eta\in\tau_\alpha\bsl(B_0\cup\dots\cup B_n)$ 

are incomparable with everything in 

$\Union_{m\leq n+1} S_m$.  More formally, let $I_n$ be (a name for)  
the set
$$
\align
I_n=
\biggl\{\eta\in{}^{<\omega}\omega&:\eta\in\tau_\alpha\bsl(B_0\cup\dots 
\cup B_n)
\ \text{and}\\
&\eta\ \text{incomparable with each}\ \eta'\in\Union_{m\leq n} S_m
\biggr\}.
\endalign
$$
The inductive construction will work with the assumption that 

$p_n\forces I_n\ \text{infinite}$ and will build $p_{n+1}$ and  
$S_{n+1}$ so that
$p_{n+1}\forces I_{n+1}\ \text{infinite}$.  Notice that both  
$S_{n+1}$ and $B_{n+1}$
are involved in the definition of $I_{n+1}$ --- so all intersections  
of
$A_\alpha=\Union_n S_n$ with $B_n$ will happen by stage $n$, and  
$A_\alpha$ must then
be almost disjoint from all of the $B_n$'s.  It is easy to check that 

$p_0\forces I_0\ \text{infinite}$, since $S_0=\emptyset$.

At stage $n+1$, construct $p_{n+1}$ and $S_{n+1}$ as follows:
Enumerate $split_{n+1}(p_n)$ as $\{\vec\sigma_i:i\leq k\}$.
Let $r'_0=p_n\restrict \vec\sigma_0$.
Extend $r'_0$ to $r''_0\leq r'_0$ so that $r''_0\forces a_0\subseteq  
I_n$ for some 

finite antichain $a_0$ satisfying $|a_0|=2^k$ and 

$a_0\cap B_{n+1}=\emptyset$ 

(this is possible because of ($*$)).
Form $r_{0,0}\leq p_n$ as follows:
$$
r_{0,0}(i)=\cases
r''_0(i)\cup \Union\{p_n(i)\restrict\sigma:
\sigma\in split_{n-i+1}(p_n(i)), \sigma\neq \vec\sigma_0(i)\},
&i\leq n+1;\\
r''_0(i)&\text{otherwise}.
\endcases
$$
So $r_{0,0}$ amalgamates $r''_0(i)$ with the ``other branches" 

of $p_n(i)$ for each $i$.
Let $r'_1=r_{0,0}\restrict \vec\sigma_1$.
The problem now is that $r'_1$ may not force the same information  
about $a_0$ that
$r''_0$ does --- in particular, $r'_1$ may have an extension which
forces all but finitely much of 

$I_n$ to be
comparable with a single node from $a_0$, ruining efforts to get 

$p_{n+1}\forces I_{n+1}\ \text{infinite}$.

So group the elements of $a_0$ into pairs
--- since the elements $t_1$ and $t_2$ in any one such pair are  
incomparable, 

$r'_1$ forces that one of $t_1$ and $t_2$ 

must have infinitely many nodes in $I_n$ incomparable to it
(if infinitely many nodes are comparable to one, then 

infinitely many of those same nodes are 

incomparable with the other).  So extend $r'_1$ to decide which of  
$t_1$ and
$t_2$ keeps infinitely much of $I_n$ incomparable to it.  Do this for
each of the $2^{k-1}$ many pairs, continuing to extend in a  
descending sequence,
so that in the end a new condition $r''_1\leq r'_1$ has 

``chosen" one node from each pair in this way.  

Form $r_{0,1}$ by amalgamating $r_{0,0}$ and $r''_1$ as follows: 

$$
r_{0,1}(i)=\cases
r''_1(i)\cup \Union\{r_{0,0}(i)\restrict\sigma:
\sigma\in split_{n-i+1}(p_n(i)), \sigma\neq \vec\sigma_1(i)\},
&i\leq n+1;\\
r''_1(i)&\text{otherwise}.
\endcases
$$
(Note that $\sigma$ may not be a splitting point of $r_{0,0}$ ---  
these are used
to preserve the {\it previous} splitting level.)
Group the remaining $2^{k-1}$ many
remaining elements of $a_0$ into pairs, and extend $r_{0,1}$ to  
``choose" between
elements of each pair, arriving at $r''_2\leq r_{0,1}$, and  
amalgamate again to get
$r_{0,2}$.
After constructing $r_{0,0}\geq r_{0,1}\geq\dots\geq r_{0,k}$, 

there will be only one element $s_0$ of
$a_0$ remaining, but the condition $r_{0,k}$ 

will force
that infinitely many elements of $I_n\bsl B_{n+1}$ 

are incomparable to $s_0$.
Now $s_0$ can be included in $S_{n+1}$ while preserving the induction  
hypothesis.

Extend $r_{0,k}\restrict \vec\sigma_1$ to decide another finite 

set $a_1$ so that  $a_1\cup\{s_0\}$ is an antichain, $a_1$ has size  
$2^k$,
$a_1\cap B_{n+1}=\emptyset$, and this extension forces  
$a_1\cup\{s_0\}\subseteq I_n$;
amalgamate the result back into $r_{0,k}$ to form $r_{1,1}$.
Group $a_1$ into pairs and extend $r_{1,1}\restrict \vec\sigma_0$ to  
decide between the
elements of the pairs
as above, amalgamating the result back into $r_{1,1}$ to obtain  
$r_{1,0}$.
Repeat as above to form $r_{1,2}$, $r_{1,3}$, etc.
Notice that the last condition $r_{1,k}$ leaves a single element
$s_1$ from the set $a_1$, and 

the condition $r_{1,k}$ force that infinitely many nodes
in $I_n\bsl B_{n+1}$ are incomparable with
{\it both} $s_0$ and $s_1$, and both can be included in $S_{n+1}$.

In general, for each $i\leq k$ form $r_{i,i}$ first to decide a  
finite set $a_i$
of size $2^k$ so that $a_i\cup\{s_0,\dots,s_{i-1}\}$ forms an  
antichain,
then form $r_{i,j}$ for $j\neq i$ in a descending sequence below
$r_{i,i}$ by grouping the set $a_i$ into pairs and extending to  
decide between 

the elements of the pairs until arriving at a single element $s_i$,
with the the $r_{i,j}$ conditions forcing now that infinitely many  
nodes
in $I_n\bsl B_{n+1}$ are incomparable with
{\it all} of the nodes $s_0,s_1,\dots,s_i$.  

At the very bottom of this finite (length $(k+1)^2$)
 descending sequence of conditions $r_{i,j}$, 

the last condition will be $r_{k,k-1}$ ---
let $p_{n+1}=r_{k,k-1}$ and $S_{n+1}=\{s_i:i\leq k\}$.

It is now easy to check that $p_{n+1}$ and $S_{n+1}$ are as desired  
--- in particular,
$p_{n+1}\forces I_{n+1}\ \text{infinite}$.
Letting  $A_\alpha=\Union_n S_n$, it is also easy to check that
$A_\alpha$ is almost disjoint from all previous $A_\beta$ (for  
$\beta<\alpha$) --- this
is what the enumeration $\{B_i:i<\omega\}$ accomplishes.
Finally, letting $q$ be the coordinatewise fusion of the sequence 

$\seq{p_n:n<\omega}$ gives $q\forces|\tau_\alpha\cap  
A_\alpha|=\omega$
since for each $n$ $q\leq p_n$ and by construction 

$p_n\forces\tau_\alpha\cap S_n\neq\emptyset$.
  This is exactly the $q$ needed for the
density argument.

Now $\Cal A=\{A_\alpha:\alpha<\omega_1\}$ is an antichain family in  
$V$.
So suppose $G$ is $\Bbb P_\omega$-generic over $V$, $\rho\in V[G]$
is an infinite antichain, and $\dot\rho$ is a name
for $\rho$.  Then $D_{\dot \rho}$ is dense, so there is a $q\in G$ 

and an $\alpha<\omega_1$ such that
$q\forces\dot \rho\cap A_\alpha\ \text{infinite}$.  Since $q\in G$,
$\rho\cap A_\alpha$ is infinite in $V[G]$. But $\rho$ was arbitrary,
so $\Cal A$ is maximal in $V[G]$.
\qed\enddemo

Recall that this Theorem deals with antichain families and $\bar\frak  
o$,
giving a stronger result than if it had dealt with off-branch  
families
and $\frak o$.

\proclaim{\Subsecnum Corollary}  

$Con(\omega_1=\frak b=\frak a=\frak d=\bar\frak o<2^\omega=\kappa)$.
\endproclaim

\demo{Proof}
Start with a model $V\satisfies ZFC+CH$.  Force to add $\kappa$-many  
Sacks reals. 

By the above Theorem, together with previous results concerning the  
effects of Sacks
forcing on the other cardinal invariants \cite{vD, BS}, 

this will give the desired model.
\qed\enddemo

\subhead\secnum Open Questions
\endsubhead

The most important of the questions seems to be the following:

\demo{\Subsecnum Question}  $Con(\frak o<\frak d)$?
\enddemo

However, the following also remains open:

\demo{\Subsecnum Question}  $Con(\frak o<\bar\frak o)$?
\enddemo

I would like to suggest a line of inquiry for answering both of these  
questions:  

Shelah's model of $Con(\frak b<\frak s)$ (see \cite{Sh}) 

gives a context in which $\frak b=\frak a=\omega_1$ and 

$\frak d={\frak a_s}=\bar\frak o= 2^\omega=\omega_2$, obtained by  
countable support 

iterated proper forcing.
(The invariant
$\frak s$ is the ``splitting number" --- see \cite{vD} for a  
definition
and basic results, \cite {BS} for the proof that $\frak s\leq\frak  
a_s$.)
If $\frak o=\omega_1$ in this construction, or there is some  
modification which will
guarantee this, then both questions are answered in the affirmative.   
I conjecture
that this technique will work.

Another avenue of inquiry open at this writing is the effects of  
various
other standard forcing notions on off-branch families (the referee
suggested Miller's superperfect forcing as an example).  

While such inquiries would prove interesting in their own right, 

they are not likely to directly address either of the above  
questions.

\enddocument